\documentclass[12pt]{amsart}
\usepackage{amssymb}
\newcommand{\A}{\mathbb{A}}
\newcommand{\Z}{\mathbb{Z}}
\newcommand{\g}{\mathfrak{g}}
\newcommand\iso{\cong}

\newcommand\Gm{\mathbb{G}_m} 
\newcommand\Ga{\mathbb{G}_a} 
\DeclareMathOperator{\Cl}{{\rm Cl}}
\newcommand{\GL}{{\rm GL}}
\newcommand{\PGL}{{\rm PGL}}
\newcommand{\dimension}{\dim}

\newtheorem{theorem}{Theorem}
\newtheorem{lemma}[theorem]{Lemma}
\newtheorem{proposition}[theorem]{Proposition}

\def\mathllap#1{\mathchoice
{\llap{$\displaystyle #1$}}%
{\llap{$\textstyle #1$}}%
{\llap{$\scriptstyle #1$}}%
{\llap{$\scriptscriptstyle #1$}}}
\def\set#1#2{\left\{\,#1\mathllap{\phantom{#2}}\mathrel{\relax}\right|\left.#2\mathllap{\phantom{#1}}\,\right\}}

\begin{document}
\title{Linear Algebraic Groups without the Normalizer Theorem}
\author{Daniel Allcock}
\address{Department of Mathematics\\University of Texas, Austin}
\email{allcock@math.utexas.edu}
\urladdr{http://www.math.utexas.edu/\textasciitilde allcock}
\thanks{Partly supported by NSF grants DMS-024512 and DMS-0600112.}
\date{August 14, 2007}

\begin{abstract}
One can develop the basic structure theory of linear algebraic groups
(the root system, Bruhat decomposition, etc.) in a way that bypasses
several major steps in the standard development, including the
self-normalizing property of Borel subgroups.
\end{abstract}

\maketitle

\noindent
An awkwardness of the theory of linear algebraic groups is that one
must develop a lot of material about general linear algebraic groups
before one can really get started.  Our goal here is to show how to
develop the root system, etc., using only the completeness of the flag
variety and some facts about solvable groups.  In particular, one can
skip over the usual analysis of Cartan subgroups, the fact that $G$ is the
union of its Borel subgroups, the connectedness of torus centralizers,
and the normalizer theorem (a Borel subgroup is
self-normalizing).  The main idea is a new approach to the
structure of rank~1 groups; the key step is lemma~\ref{lem-rank-1-implies-exactly-2-Borels}.

\smallskip
All algebraic geometry is over a fixed algebraically closed field.
$G$ always denotes a connected linear algebraic group with Lie algebra
$\g$, $T$ a maximal torus, and $B$ a Borel subgroup containing it.  We
assume the structure theory for connected solvable groups, and the
completeness of the flag variety $G/B$ and some of its consequences.
Namely: that all Borel subgroups (resp. maximal tori) are conjugate;
that $G$ is nilpotent if one of its Borel subgroups is nilpotent; that
$C_G(T)_0$ lies in every Borel subgroup containing $T$; and that
$N_G(B)$ contains $B$ of finite index and (therefore) is
self-normalizing.  We also assume known that the centralizer of a
torus has the expected dimension, namely, that of the subspace of $\g$
where the torus acts trivially.  For these results we refer to
Borel \cite{borel}, Humphreys \cite{humphreys} and Springer \cite{springer}.

In section~\ref{sec-solvable} we develop a few properties of solvable groups, and
in section~\ref{sec-rank-1} we treat the structure of rank~1 groups.  The root
system, etc., can then be developed in essentially the standard way,
so after the rank~1 analysis we restrict ourselves to brief comments.
We are grateful to J.~Humphreys, G.~McNinch and T.~Springer for their
helpful comments. 

\section{Lemmas about solvable groups}
\label{sec-solvable}

\noindent
First we recall from \cite[\S13.4]{springer} the groups that we call the
positive and negative subgroups of $G$.
Fix a 1-parameter group $\phi:\Gm\to G$.  For $g\in G$  and
$\lambda\in\Gm$, define
$\Cl_g(\lambda)=\phi(\lambda)g\phi(\lambda)^{-1}$ and 
$$
G^+=\set{g\in G}{\lim_{\lambda\to0}\Cl_g(\lambda)=1}.
$$ 
Of course, the condition $\lim_{\lambda\to0}\Cl_g(\lambda)=1$ means
that $\Cl_g:\Gm\to G$ extends to a regular map $\Gm\cup\{0\}\to G$
sending $0$ to~$1$.
We call $G^+$ the positive part of $G$ (with respect to
$\phi$).  It is a group because
$$
\lim_{\lambda\to0}\phi(\lambda)gg'\phi(\lambda)^{-1}
=
\lim_{\lambda\to0}\phi(\lambda)g\phi(\lambda)^{-1}
\cdot\lim_{\lambda\to0}\phi(\lambda)g'\phi(\lambda)^{-1}
$$
and
$$
\lim_{\lambda\to0}\phi(\lambda)g^{-1}\phi(\lambda)^{-1}
=
\Bigl(\lim_{\lambda\to0}\phi(\lambda)g\phi(\lambda)^{-1}\Bigr)^{-1}
$$ when the limits on the right hand sides exist.  It is closed and
connected because it is generated by irreducible curves
containing~$1$.  Of course, there is a corresponding subgroup $G^-$
called the {negative part of $G$}, got by considering limits as
$\lambda$ approaches~$\infty$.  All of our discussion applies equally
well to $G^-$.
The key properties of $G^\pm$ are that they are unipotent
and ``large'':

\begin{proposition}
\label{prop-pos-and-neg-groups}
$G^+$ is unipotent, and every weight of $\Gm$ on the lie algebra of
$G^+$ is positive.  If $G$ is solvable, then the lie algebra of $G^+$
contains all the positive weight spaces for $\Gm$ on $\g$.
\end{proposition}

\noindent
We remark that the solvability hypothesis in the last part is
unnecessary, but we will need only the solvable case.  The general
case requires the structure theory which depends on theorem~\ref{thm-rk-1-implies-PGL2}.

\begin{proof}
See \cite[theorem 13.4.2]{springer} for the first claim; the idea is
to embed $G$ into $\GL_n$ and diagonalize the $\Gm$ subgroup.  We
will prove the second claim, since the proof in
\cite{springer} includes the nonsolvable case and therefore relies on
the Bruhat decomposition.  We will use induction on $\dimension
G_u$;  the $0$-dimensional case is trivial.  

So suppose $\dimension
G_u>0$ and choose a connected subgroup $N$ of $G_u$, normal in $G$,
with $G_u/N\iso\Ga$.  We write $\pi$ for the associated map
$G_u\to\A^1$.  By the action of  $\Gm$  on $G_u/N$, we have $\pi\circ\Cl_g(\lambda)=\pi(g)\cdot \lambda^n$ for some $n\in\Z$.
If $n\leq0$ then $\pi\circ\Cl_g$  does not extend to a regular map $\Gm\cup\{0\}\to\A^1$ sending
$0$ to $0$ unless $\pi(g)=0$, so a group element outside $N$ cannot
lie in $G^+$.  So $G^+=N^+$ and we use induction.  

So suppose $n>0$, i.e., the action is by a positive character.
Choose a linear representation $V$ of $\Gm$ and an embedding (as a
variety) of $G_u$ into $PV$ which is equivariant with respect to the
conjugation action of $\Gm$.  Choose a $\Gm$-invariant linear subspace
of $PV$ containing~$1\in G_u$, whose
tangent space there is complementary to that of $N$.  By dimension
considerations, its intersection with $G_u$ contains  an invariant
curve $C$ that passes through $1$ and does not lie in $N$.  
By passing to a component we may assume
that $C$ is irreducible, so it is the closure of the orbit of some
$g\in C$.  The map $\Cl_g:\Gm\to G_u$ extends to a regular map from
$P^1$ to $PV$.  Because $1$ lies in the
closure of the orbit, $\Cl_g$ sends at least one of $0$ or $\infty$ to
$1$.  It cannot send $\infty$ to $1$, because $\pi\circ\Cl_g(\lambda)=
\hbox{(nonzero constant)}\cdot\lambda^n$ for some $n>0$, which admits
no regular extension $\Gm\cup\{\infty\}\to\A^1$.  Therefore
$\lim_{\lambda\to0}\Cl_g(\lambda)=1$, so $g\in G^+$.  This shows that
$G^+$ projects onto $G_u/N$.  It also contains $N^+$, to which the
inductive hypothesis applies.  The proposition follows.
\end{proof}

An immediate consequence is that a connected solvable group $G$ is generated by
its subgroups $G^+$, $G^-$ and $C_G(\phi(\Gm))_0$, since together their
Lie algebras span $\g$.
Next, we need a theorem on orbits of solvable groups.  The
proof does not use the structure theorem for solvable
groups, and indeed can be used to prove it.   It can
be simplified slightly if one assumes the structure theorem.

\begin{theorem}
\label{thm-solvable-homogeneous-spaces-lack-complete-subsets}
If $G$ is solvable and acts on a variety, then no orbit contains any
complete subvariety of dimension${}>0$.
\end{theorem}

\begin{proof}
We assume the result known for solvable groups of smaller dimension
than $G$.  (The $0$-dimensional case is trivial.)  If $X$ is the
variety, $x\in X$, $G_x$ its stabilizer and $Y$ its orbit, then the
natural map $G/G_x\to Y$ is generically finite (since $G/G_x$ and $Y$
have the same dimension), hence finite (by homogeneity).  
If $Y$ contained a complete subvariety of positive dimension then the
preimage in $G/G_x$ would also be complete.  So it suffices to treat
the case $X=G/G_x$.

We consider three cases. First, if $G_x$ contains $[G,G]$ then it is
normal, so $G/G_x$ is an affine variety and cannot contain a complete
subvariety of positive dimension.  Second, if $G_x$ surjects to
$G/[G,G]$, then $[G,G]$ acts transitively on $G/G_x$, and by the
inductive hypothesis applied to $[G,G]$, $G/G_x$ cannot contain a
complete subvariety of positive dimension.  Finally, suppose $G_x$
neither contains $[G,G]$ nor surjects to $G/[G,G]$.  Set $H$ equal to
the group generated by $G_x$ and $[G,G]$.  We will use the fact that
$G/G_x$ maps to $G/H$ with fibers that are copies of $H/G_x$.  As in
the first case, $G/H$ is an affine variety, so any complete subvariety
of $G/G_x$ lies in the union of finitely many copies of $H/G_x$.  But
the inductive hypothesis applied to $H$ shows that every complete
subvariety of $H/G_x$ is a finite set of points.  Therefore the same
conclusion applies to $G/G_x$.
\end{proof}

\section{Rank One Groups}
\label{sec-rank-1}

\noindent
In this section, $G$ is connected and non-solvable of rank~1.  To goal is:

\begin{theorem}
\label{thm-rk-1-implies-PGL2}
$G$ modulo its unipotent radical admits an isogeny to $\PGL_2$.
\end{theorem}

There is a standard argument that reduces this to proving that $T$
lies in exactly two Borel subgroups.  We must modify this slightly
because we are not assuming the normalizer theorem. 
Since $N:=N_G(B)$ is
self-normalizing, it fixes only one point of $G/N$, so
the stabilizers of distinct points of $G/N$ are the normalizers of distinct
Borel subgroups.  The fixed points of $T$ in
$G/N$ correspond to Borel subgroups that $T$ normalizes, hence lies in.
Now we use the theorem that a torus acting on a $d$-dimensional
projective variety has at least $d+1$ fixed points.  Since $G$ is not solvable,
$G/B$ has dimension${}>0$, so $G/N$ does too.  Therefore $T$ lies in
at least two Borel subgroups.  And if we prove that it lies in exactly
two, then we can also deduce $\dimension G/N=1$.  Then it is easy to see that $G/N\iso P^1$ and
prove theorem~\ref{thm-rk-1-implies-PGL2}.  
So our aim is to prove that $T$ lies in exactly two Borel
subgroups.  

\medskip
Using the positive and negative subgroups, we will construct two
Borel subgroups containing $T$, and   then show that there are no
more.  Suppose $\phi:\Gm\to T$ is a parametrization of $T$ (meaning
$\phi$ is an isomorphism) and $B$ a Borel subgroup containing $T$.  Call $B$
{positive} (with respect to $\phi$) if it contains $G^+$ and
{negative} if it contains $G^-$.  Obviously, $B$ is positive with
respect to one parametrization of $T$ if and only if it is negative
with respect to the other.  Here are the basic properties of positive
and negative Borel subgroups.

\begin{lemma}
\label{lem-pos-and-neg-Borels-properties}
Suppose  $\phi:\Gm\to T$ a parametrization of the maximal torus $T$.  Then
\begin{enumerate}
\item $T$ lies in a positive Borel subgroup and in a negative Borel subgroup;
\label{item-pos-neg-Borels-exist}
\item if $B$ (resp. $B'$) is a positive (resp. negative) Borel
  subgroup containing $T$, then every Borel subgroup containing $T$
  lies in $\langle B,B'\rangle$;
\label{item-two-borels-contain-all}
\item every Borel subgroup containing $T$ is either positive or negative, but
  not both; and
\label{item-Borel-pos-or-neg-not-both}
\item $N_G(T)$ contains an element acting on $T$ by inversion.
\label{item-Weyl-group-contains-negation}
\end{enumerate}
\end{lemma}

\begin{proof}
\eqref{item-pos-neg-Borels-exist}  $G^+$ is connected, 
unipotent and normalized by $T$.  Therefore $TG^+$ lies in some
Borel subgroup, which is then positive.  And similarly for $G^-$.

\eqref{item-two-borels-contain-all}  Suppose $B''$ is a Borel subgroup
containing $T$.  Then ${B''}^+$ lies in $B$ since $B$ is positive,
${B''}^-$ lies in $B'$ since $B'$ is negative, and $C_{B''}(T)_0$ lies
in both $B$ and $B'$ because $C_G(T)_0$ lies in every Borel subgroup
containing $T$.  Since $B''$ is generated by ${B''}^+$, ${B''}^-$ and $C_{B''}(T)_0$, 
it lies in $\langle B,B'\rangle$.  

\eqref{item-Borel-pos-or-neg-not-both} The condition that a Borel subgroup
containing $T$ contains at least one of $G^\pm$ is independent under
conjugation by $N(T)$.  Since $N(T)$ acts transitively on the
Borel subgroups containing $T$, and we have exhibited positive and negative
Borel subgroups, every Borel subgroup containing $T$ contains at least one of $G^\pm$
and hence is either positive or negative.  If a Borel subgroup $B$
containing $T$ were both positive and negative, then \eqref{item-two-borels-contain-all} would imply
that it is the only Borel subgroup containing $T$, contradicting the
fact that $T$ lies in at least~2 Borel subgroups.

\eqref{item-Weyl-group-contains-negation} By \eqref{item-pos-neg-Borels-exist}, positive and
negative Borel subgroups exist, and by \eqref{item-Borel-pos-or-neg-not-both} they are distinct.
The result follows from the fact that $N_G(T)$ acts transitively on
the Borel subgroups containing $T$.
\end{proof}

Now we can give the key step in our approach to the structure theorem
for rank~1 groups.

\begin{lemma}
\label{lem-rank-1-implies-exactly-2-Borels}
Every
maximal torus of $G$ lies in exactly two Borel subgroups, one positive and one
negative. 
\end{lemma}

\begin{proof}
By induction on the dimension of a Borel subgroup.  If it has
dimension~1 then it is abelian, so $G=B$.  Therefore the
base case is dimension~2.  Choose a maximal torus $T$ and a
parametrization of it.  We already know that $T$ lies in a positive
and a negative Borel subgroup.  The key point is that any two positive Borel subgroups
coincide.  For otherwise their unipotent radicals would be distinct
subgroups of $G^+$, hence generate a unipotent group of
dimension${}>1$.  This is impossible because $\dimension B_u=1$.
Similarly, there is only one negative Borel subgroup and the base case is
proven.

Now we prove the inductive step; suppose $B$ has dimension at
least $3$, and suppose without loss that it is positive.  Consider the action of $B$ on $G/N$; there is a unique
fixed point because the only Borel subgroup that $B$ normalizes is itself.
Consider an orbit of minimal positive dimension.  By
theorem~\ref{thm-solvable-homogeneous-spaces-lack-complete-subsets},
it contains no complete subvarieties of dimension${}>0$.  On the other
hand, its closure is complete and is got by adjoining
lower-dimensional orbits.  By minimality, this means that its closure
is got by adjoining a single point, so the orbit is a curve.  Therefore there exists a Borel subgroup $B'$ for
which $B\cap N(B')$ has codimension~1 in $B$.  That is,
$I:=(B\cap B')_0$ has codimension~1 in each of $B$ and $B'$.  There are two
possibilities: $I=B_u=B'_u$, or $I$ contains a torus.
In the first case, $\langle B,B'\rangle$ normalizes $I$, and a
Borel subgroup in $\langle B,B'\rangle/I$ has no unipotent part.  This
forces $\langle B,B'\rangle$ to be solvable, which is impossible.

Therefore $I$ contains a torus $T$, and $I_u$ has codimension~1 in
each of $B_u$ and $B_u'$.  Now, $T$ normalizes $B_u$ and $B_u'$, hence
their intersection, hence $I_u$.  Also, $B_u$ normalizes $I_u$ because
it is only one dimension larger and is nilpotent.  Similarly for
$B_u'$.  Therefore $\langle B,B'\rangle$ normalizes $I_u$,
which has dimension${}>0$ since $\dimension B>2$.  We apply
induction to $\langle B,B'\rangle/I_u$, and then pull back to conclude
the following.  $B$ and $B'$ are the only Borel subgroups of $\langle B,B'\rangle$
containing $T$, and they are exchanged by an element of $N_{\langle
  B,B'\rangle}(T)$ that inverts $T$.  This implies that $B'$ is a
negative Borel subgroup of $G$.
Finally, lemma~\ref{lem-pos-and-neg-Borels-properties}\eqref{item-two-borels-contain-all}
implies that any Borel subgroup of $G$ containing $T$ lies in $\langle
B,B'\rangle$, hence equals $B$ or $B'$.
\end{proof}

This lemma implies theorem~\ref{thm-rk-1-implies-PGL2}, and from then on one can follow
the standard development.  We make only the following remarks.


{\it Bruhat decomposition:\/} in the absence of the normalizer
theorem, one should define the Weyl group $W$ as the subgroup of
$N_G(T)/C_G(T)$ generated by the reflections coming from roots.  Then
one can prove $G=BWB$ as in \cite[\S14]{borel}, \cite[\S28]{humphreys} or
\cite[\S8.3]{springer}. 

{\it Normalizers:\/} The theorem $N_G(B)=B$ follows immediately from the Bruhat
decomposition, and implies that $W$, as defined here, is all of
$N_G(T)/C_G(T)$, so that our definition agrees with the usual one.

{\it Connectedness of torus centralizers:\/} this can be deduced from
the Bruhat decomposition and a standard fact about  reflection
groups: the pointwise stabilizer of a linear subspace is generated by
the reflections that stabilize it.

\end{document}